\documentclass[12pt]{article}
\usepackage{amsfonts}
\usepackage{graphicx}
\usepackage{amsmath}
\usepackage{fancyhdr}
\usepackage{amsmath}
\usepackage{amssymb}
\usepackage{amsthm}
\usepackage{amsbsy}
\usepackage{makeidx}
\usepackage[cp866]{inputenc}
\begin{document}

\begin{center}
\textbf{ Topp-Leone generated q-exponential distribution and its applications
}

\bigskip

Nicy Sebastian, Rasin R. S. and Silviya P. O.

\vspace{3mm}  Department of Statistics, St.Thomas College, Thrissur,
Kerala, India-680001

e-mail: nicycms@gmail.com, rasinrs@gmail.com

\bigskip
\end{center}

\begin{abstract}

 Topp-Leone distribution is a continuous model distribution used for modelling lifetime phenomena. The main purpose of this paper is to introduce a new framework for generating
lifetime distributions, called the Topp-Leone generated q-exponential family of distributions. Parameter estimation using maximum likelihood method and simulation results to assess effectiveness of the distribution are discussed. Different informative and non-informative priors are used to estimate the shape parameter of q extended Topp-Leone generated exponential distribution under normal approximation technique. We prove empirically the importance and flexibility of the new model in model building by using a real data set.

\vskip.2cm
 \noindent\textit{Keywords:} Topp-Leone distribution, beta-generated,
generalized exponential,  parameter estimation, simulation. \\

 \noindent\textit{MSC (2010)} 60E05, 33B20, 62G05, 62N05,  68U20.

\end{abstract}

\section{Introduction}
\vskip.3cm
\noindent
	There are many statistical distributions which plays an important role in modeling survival and life time data such as exponential, weibull, logistic etc. Almost all these  distributions with unbounded support. But there are situations in real life, in which observations can take values only in a limited range such as percentages, proportions or fractions. Papke and Wooldridge (1996) claims that in many economic settings, such as fraction of total weekly hours spent working, pension plan participation rates, industry market shares, fraction  of land area allocated to agriculture etc., the variable bounded between zero and one. Thus it is important to have models defined on the unit interval in order to have reasonable results. Also different authors refer to continuous models with finite support in order to describe  life time data, in reliability analysis.
 It is well known that beta  distribution is the most used distribution to model continuous variables in the unit interval. This distribution is popular in the field of engineering, economics, biology, ecology etc. due to the great flexibility of its density function. But due to the fact that  its distribution function cannot be expressed in closed form and it involves the incomplete beta function, the mathematical formulation is found to be difficult. However, several authors have proposed alternatives to the beta distribution by recovering the distribution proposed by Kumaraswamy in 1980.
\vskip.2cm
\noindent

A new distribution was introduced in 1955, called Topp Leone (TL) distribution, defined on finite support, proposed  Topp and Leone and used it as a model for failure data. A random variable $X$ is distributed as the TL with parameter $\alpha$ denoted by $x\sim TL _ {(\alpha)}$, with a cumulative distribution function
\begin{equation}\label{eq 1}
F_{TL}(x)=x^\alpha(2-x)^\alpha ,0<x<1,\alpha >0.
\end{equation}
The corresponding probability function is
\begin{equation}\label{eq 2}
f_{TL}(x)=2\alpha x^{\alpha-1}(1-x)(2-x)^{\alpha-1}	.
\end{equation}

  Topp Leone distribution provides closed forms of cumulative density function (cdf) and the probability density function (pdf) and describes empirical data with J-shaped histogram such as powered tool band failures, automatic calculating machine failure.  The Topp Leone distribution had been received little attention until Nadarajah and Kotz (2003) discovered it. They studied about some properties of TL distribution and provided its moments, central moments and characteristic function. Ghitany et al. (2005) provided some reliability measures of TL distribution such as a hazard function, mean residual life, reversed hazard rate, expected inactivity time, an its stochastic orderings. A discussion on kurtosis of the TL distribution was reported by Kotz and Seier (2002).

\noindent Lifetime data plays an important role in a wide range of applications such as medical, engineering an social sciences.  When there is a need for more flexible distributions, almost all researchers are about to use the new one with more generalization. An excellent review of Lee et al.(2013) has provided through knowledge of several methods for generating families of continuous univariate distributions. They classify these methods by years before and after 1980. According to their work, there are some general methods introduced prior to 1980, and they may be summarized as the method of differential equation, method of transformation (also known as translation), method of quantile. Since 1980, methodologies of generating new distribution shifted to adding new parameters to an existing distribution. According to Lee et al.(2013) some noticeable developments after 1980, are method of generating skew distributions, beta generated method, method of adding parameters, transformed-transformer method, and composite method.

\noindent The beta generated (BG) family of distributions belongs to a parameter adding mathod (Lee et al., (2013)). For an arbitrary distribution with a cumulative distribution function (cdf) $G(\cdot)$ and a probability density function (pdf) $g(\cdot)$, this method generates it by letting $x=G^{-1}(B)$ where $B$ is the beta function, $B \sim Beta (a,b)$ (see Alexander et al., 2012). Some existing distributions incorporated with BG family will have two additional parameters, which are the parameters of beta distribution.
The cdf of beta generated random variable $X$ is defined as
\begin{equation*}
F_{BG}(x)=\int_{0}^{G(x)} h(t)dt
\end{equation*}
where $h(t)$ is the pdf of beta random variable and $G(x)$ is the cdf of any arbitrary random variable. Thus the cdf of beta generated random variable is
\begin{equation} \label{eq 3}
F_{BG}(x)=\frac{1}{B(a,b)}\int_{0}^{G(x)}t^{a-1}(1-t)^{b-1}dt, a>0, b>0.
\end{equation}
  The pdf corresponding to \ref{eq 3}) is given by
\begin{equation} \label{eq 4}
f_{BG}(x)=\frac{1}{B(a,b)} g(x) G(x)^{a-1} (1-G(x))^{b-1},
\end{equation}
where $ B(a,b)=\frac{\Gamma a \Gamma b}{\Gamma (a+b)} $ is the beta function and $ \Gamma (\cdot) $ is the gamma function.

\noindent
Instead of using beta distribution as generator, we use TL distribution,  as generator distribution, and we obtain TLG family of distribution.
Then relation of a random variable $X $having the TLG distribution and a random variable $T$ having TL distribution is $X=G^{-1}(T),$ with $ T \sim  TL(\alpha)$. This relation demonstrates that the pdf of TL distribution, (\ref{eq 2}), is transformed into a new pdf through the function $G({^.})$
\begin{equation} \label{eq 6}F_{TLG}(x)=2\alpha \int_{0}^{G(x)}{t^{\alpha-1}(1-t)(2-t)^{\alpha-1}dt}
	=G(x)^\alpha(2-G(x))^\alpha.
	\end{equation}
By differentiating, we get the corresponding pdf,
\begin{equation} \label{eq 7}
f_{TLG}(x)=2\alpha g(x)(1-G(x))G(x)^{\alpha-1}(2-G(x))^{\alpha-1},\alpha>0.
\end{equation}
\begin{figure}[h!]
	\begin{center}
		~~~~~ \resizebox{7cm}{7cm}{\includegraphics{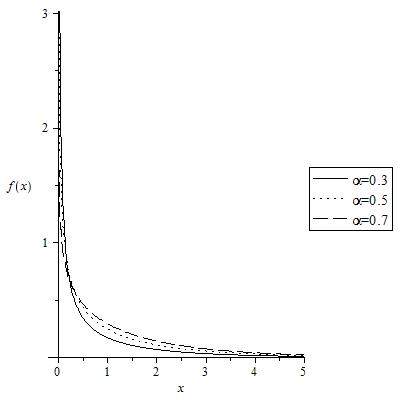}}
		~~~~ \resizebox{7cm}{7cm}{\includegraphics{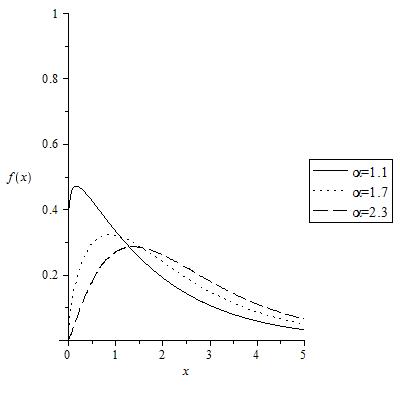}}\\
		\caption{\small  Plots of density function of TLE distribution with $\lambda=0.3$    \label{ppplot2}}
	\end{center}
\end{figure}
\noindent In reliability analysis, a frequently used distribution is exponential distribution having the  characterizing property of constant hazard function. Due to this, exponential distribution is sometimes not suitable for analyzing data. This implies the need for more generalization. In such situations we use distribution called Topp-Leone Exponential distribution (TLE). TLE distribution comes as the combination of TL distribution and exponential distribution. Here TL distribution is the generator and exponential is the parent distribution.
For creating TLE distribution, we need cdf G(x) and pdf g(x) of exponential distribution, \begin{equation}\label{eq 10}
G(x)=1-{\rm e}^{-\lambda x}; x>0,\lambda >0,
\end{equation}
and
\begin{equation} \label{eq 11}
g(x)=\lambda {\rm e}^{-\lambda x}.
\end{equation}
The TLE distribution is obtained by taking (\ref{eq 10}) and (\ref{eq 11}) into (\ref{eq 6}) and (\ref{eq 7}).
Sangsanit and Bodhisuwan (2016) presented the Topp-Leone generated exponential (TLE) distribution as an example of the Topp-Leone generated
distribution. A random variable $X$ possessing TLE distribution having cdf and probability function defined respectively  as
\begin{equation*} F_{TLE}(x)=(1-exp(-\lambda x))^\alpha(2-(1-exp(-\lambda x)))^\alpha =(1-exp(-2\lambda x))^\alpha
	\end{equation*}	
and
\begin{equation*}
f_{TLE}(x)=2\alpha \lambda exp(-2 \lambda x)(1- exp(-2\lambda x))^{\alpha-1}.
\end{equation*}
Where $ \alpha $ is the shape parameter and $ \lambda $ is the scale parameter.
The survival and hazard function of TLGE distribution is given as
\begin{center}
	$ S(x)=1-(1-exp(-2\lambda x))^\alpha, $
\end{center}
and
\begin{center}
	$ h(x)=\frac{2\alpha\lambda exp(-2\lambda x)(1-exp(-2\lambda x))^{\alpha-1}}{1-(1-exp(-2\lambda x))^\alpha} $ .
\end{center}
respectively. It will be useful to consider the shape of the hazard function, in reliability analysis, to select appropriate distribution since it is an important measure of aging. For TLE distribution, depending on the values of the parameters,  it can have constant, increasing, and decreasing  hazard function. Plots of probability function and hazard function are  respectively given in  Figure 1 and Figure 2.

\begin{figure}[h!]
	\begin{center}
		~~~~~ \resizebox{7cm}{7cm}{\includegraphics{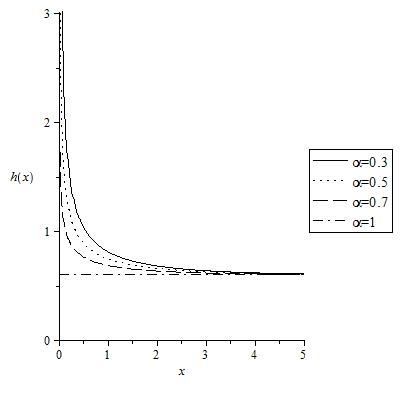}}
		~~~~ \resizebox{7cm}{7cm}{\includegraphics{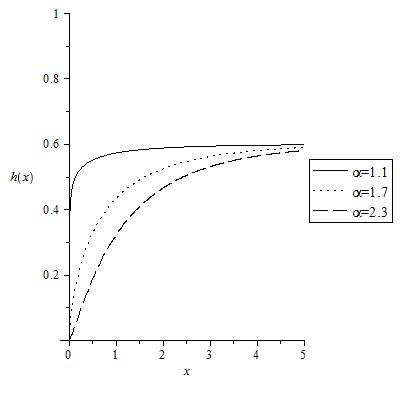}}\\
		\caption{\small  (b) $h(x)$ of TLE with $\lambda=0.3$, $\alpha \leq1$ \hskip.5cm (c) $h(x)$ of TLE with $\lambda=0.3$, $\alpha >1$ 	 \label{ppplot2}}
	\end{center}
\end{figure}

\section{q-Exponential Distribution}
\vskip.3cm
\noindent
Various entropy\index{entropy} measures have been developed by
mathematicians and physicists to describe  several phenomena,
depending on the field and the context in which it is being used.
In statistical mechanics, Maxwell-Boltzmann distribution maximizes the
Boltzmann-Gibbs entropy under appropriate constraints (Gell-Mann and Tsallis (2004)).
Given a probability distribution $P = \{p_i\} (i = 1,\ldots,m),$ with $p_i$ representing the
probability of the system to be in the $i$th microstate, the Boltzmann-Gibbs entropy is
\begin{equation}\label{eq:0}S=-k\sum_{i=1}^mp_i \ln p_i, \end{equation}
 where $k$ is the Boltzmann constant and $m$ the total number
of microstates. If all states are equally probable it leads to the Boltzmann principle
$S = k \ln W (m = W)$. Boltzmann-Gibbs entropy is equivalent to Shannon's entropy if
$k = 1$. If we consider such a system in contact with a thermostat
then we obtain the usual Maxwell-Boltzmann distribution for the possible
states by maximizing the Boltzmann-Gibbs entropy $S$ with the normalization
and energy constraints. However, in nature many systems show distributions
which differ from the Maxwell-Boltzmann distribution. Tsallis (1988),  introduced a
  generalization of  the Boltzmann-Gibbs \index{Boltzmann-Gibbs entropy} entropy\index{entropy}.
  The $q$-entropic function is of the form
   $$ S_{T}(P)= k \frac{W^{1-q}-1}{1-q}.$$
      By  maximizing  Tsallis\index{Tsallis entropy} entropy\index{entropy}, subject to certain constraints, leads to
   the  Tsallis distribution, also known as $q$-exponential\index{$q$-exponential} distribution, which has the form $f(x) = c [1-(1-q)x]^{\frac{1}{(1-q)}}$
   where $c$ is the normalizing constant. Various applications and generalizations of the $q$-exponential\index{$q$-exponential} distribution are given in Picoli et al (2003).
    In the limit $q\rightarrow 1$, $q$-entropy converges to Boltzmann-Gibbs entropy\index{Boltzmann-Gibbs entropy}.

 An important characteristic of $q$-exponential distribution is that it has two parameters $q$ and $\lambda$ providing more flexibility with regard to its decay, differently from exponential distribution.
The $q$ exponential distribution is defined by its pdf and cdf as,
\begin{equation}\label{eq 6a}
f_{qE}(x)=(2-q)\lambda [1-(1-q)\lambda x]^\frac{1}{1-q}, x>0, \lambda >0, q<2, q\neq 0
\end{equation}

\begin{equation}\label{eq 7b}
F_{qE}(x)=1-[1-(1-q)\lambda x]^{\frac{(2-q)}{1-q}}.
\end{equation}
The parameter $q$ is known as entropy index.
As $q\rightarrow1$ $q$-exponential distribution becomes exponential distribution. In that sense $q$-exponential distribution is a generalization of exponential distribution. The parameters $q$ and $\lambda$ determine how quickly the pdf decays. In the reliability context, an important characteristic of the $q$-exponential distribution is its hazard rate. It is given as
\begin{eqnarray}
\nonumber
h_{qE}(x)
&=&\lambda (2-q) [1-(1-q)\lambda x]^\frac{1}{1-q}
\end{eqnarray}
For $q$-exponential distribution, the hazard rate is not necessarily constant as in exponential distribution. For $1<q<2$,	$h_{qE}(x)$ is a decreasing monotonic function, while for $q<1$, $h_{qE}(x)$ increases monotonically.

\vskip.3cm
\section{Topp-Leone q-Exponential Distribution}
\vskip.3cm

In this section  we introduce the Topp-Leone generated q- exponential(TLqE) distribution by combining the TL distribution with q-exponential distribution. Substituting (\ref{eq 6a}) and (\ref{eq 7b}) in (\ref{eq 6}) and (\ref{eq 7})respectively we will get the
 distribution function and density function of TLqE distribution as

\begin{eqnarray}
\nonumber
F_{TLqE}(x)&=&\{1-[1-(1-q)\lambda x]^{\frac{2-q}{1-q}}\}^{\alpha}\{2-\{1-[1-(1-q)\lambda x]^{\frac{2-q}{1-q}} \} \}^{\alpha} \\
\nonumber \\
&=&\{ 1-[1-(1-q) \lambda x]^{2( \frac{2-q}{1-q})}\}^{\alpha},  x>0,\lambda,\alpha>0, q<2,q\neq0 ,
\end{eqnarray}

and
\begin{eqnarray}
\nonumber
f_{TLqE}(x)&=&2 \alpha (2-q) \lambda[1-(1-q) {\lambda x}]^{\frac{1}{1-q}} [1-(1-q) \lambda x]^{2 (\frac{2-q}{1-q})}\\
\nonumber
& &\times \{ 1-[1-(1-q) \lambda x ]^{2 (\frac{2-q}{1-q})} \}^{\alpha -1}  \\
&=&2 \alpha \lambda(2-q)[1-(1-q) \lambda x]^{\frac{3-q}{1-q}}\{ 1-[1-(1-q) \lambda x]^{2 (\frac{2-q}{1-q})}   \}^{\alpha-1}.
\end{eqnarray}

\begin{figure}[h!]
	\begin{center}
		~~~~~ \resizebox{7cm}{7cm}{\includegraphics{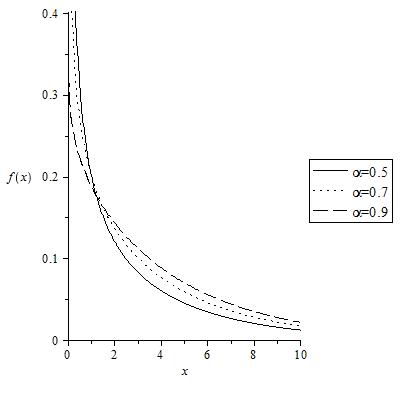}}
		~~~~ \resizebox{7cm}{7cm}{\includegraphics{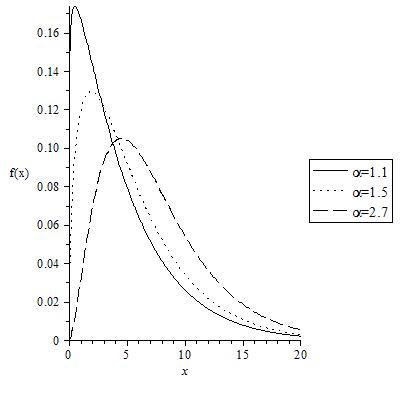}}\\
		\caption{\small Plots of f(x) of TLqE distribution for $\alpha<1$ (left) and for $\alpha>1$ (right) with $\lambda=0.1, q=0.9$   \label{ppplot2}}
	\end{center}
\end{figure}
\begin{figure}[h!]
	\begin{center}
		~~~~~ \resizebox{8cm}{8cm}{\includegraphics{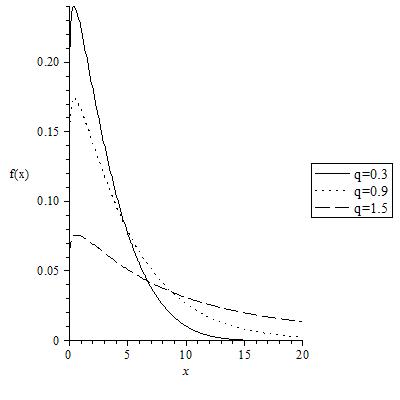}}
		\caption{\small Plots of f(x) of TLqE distribution for different values of the parameter q with $\alpha=1.1,$ $\lambda=0.1$   \label{ppplot2}}
	\end{center}
\end{figure}

In Figure 3 and Figure 4, we can see the plots of density function of TLqE for different values of the shape parameters $\alpha$ and $q$.
The survival function, the probability density function and the Hazard function are the three important functions that characterize the distribution of the survival times. Here
\begin{center}
	$ S(x)=1-\{ 1-[1-(1-q) \lambda x]^{2( \frac{2-q}{1-q})}\}^{\alpha} $,
\end{center}
and
\begin{center}
	$	h(x)= \frac{2 \alpha \lambda(2-q)[1-(1-q) \lambda x]^{\frac{3-q}{1-q}}\{ 1-[1-(1-q) \lambda x]^{2 (\frac{2-q}{1-q})}   \}^{\alpha-1}}{1-\{ 1-[1-(1-q) \lambda x]^{2( \frac{2-q}{1-q})}\}^{\alpha}} $
\end{center}
respectively are the survival and the hazard function of TLqE distribution. Figure 5, gives the plots of $h(x)$ of TLqE distribution for $\alpha\leq1$ (left) and for $\alpha>1$ (right) with $\lambda=0.1, q=0.9$. The cumulative hazard function is given as
\begin{eqnarray}
\nonumber
H(x)&=&\int_{0}^{x}h(t)dt\\
\nonumber
&=&- \ln\{ 1-\{ 1-[1-(1-q) \lambda x]^{2( \frac{2-q}{1-q})}\}^{\alpha} \}.
\end{eqnarray}

\begin{figure}[h!]
	\begin{center}
		~~~~~ \resizebox{7cm}{7cm}{\includegraphics{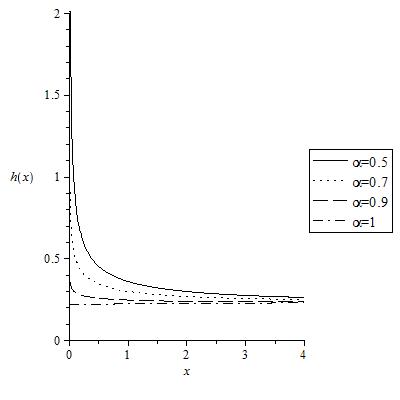}}
		~~~~ \resizebox{7cm}{7cm}{\includegraphics{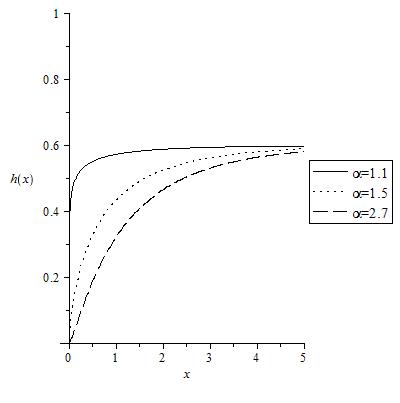}}\\
		\caption{\small Plots of h(x) of TLqE distribution for $\alpha\leq1$ (left) and for $\alpha>1$ (right) with $\lambda=0.1, q=0.9$   \label{ppplot2}}
	\end{center}
\end{figure}

\subsection{Maximum Likelihood Estimation}
Here we derive the maximum likelihood estimate of the unknown parameter vector $\theta=(\alpha,\lambda,q)^T$. Suppose $ x_{1},x_{2},...,x_{n} $ be a random sample of size $n$ taken from TLqE distribution. We have the pdf of TLqE distribution is
\begin{equation*}
f_{i}(x,\alpha, \lambda,q)=2 \alpha \lambda (2-q)[1-(q-1) \lambda x_{i}]^{\frac{3-q}{1-q}} \{ 1-[1-(1-q)\lambda x_{i}]^{\frac{2(2-q)}{(1-q)}} \}^{\alpha -1}
\end{equation*}
Then the likelihood function is written as
\begin{eqnarray}\label{eq:lf}
\nonumber
L(\theta)&=&\prod_{i=1}^{n} f_{i}(x,\alpha, \lambda,q) \\
&=&[2 \alpha \lambda (2-q)]^{n} \prod_{i=1}^{n} [1-(1-q)\lambda x_{i}]^{\frac{(3-q)}{(1-q)}} \{ 1-[1-(1-q)\lambda x_{i}]^{\frac{2(2-q)}{(1-q)}} \}^{\alpha -1}
\end{eqnarray}
and the log likelihood function is
\begin{eqnarray}
\nonumber
\ell(\theta)&=& n \ln 2 + n \ln\alpha + n \ln \lambda + n \ln (2-q) \\
\nonumber
& &\frac{(3-q)}{(1-q)} \sum_{i=1}^{n} \ln [1-(1-q)\lambda x_{i}] + (\alpha-1) \sum_{i=1}^{n} \ln \{ 1-[1-(1-q)\lambda x_{i}]^{\frac{2(2-q)}{(1-q)}} \}
\end{eqnarray}
Let $ h_{i}(x,\lambda ,q) = 1-(1-q)\lambda x_{i}$
$ r=\frac{2(2-q)}{(1-q)} $
and $ r-1= \frac{2(2-q)}{(1-q)} -1 =\frac{3-q}{1-q}$.
Thus $ \ell(\theta) $ can be written as
\begin{eqnarray}
\nonumber
\ell(\theta)&=&n \ln 2 + n \ln \alpha + n \ln \lambda + n \ln (2-q)  \\
\nonumber
& &+ (r-1) \sum_{i=1}^{n} \ln h_{i}(x,\lambda, q) + (\alpha-1) \sum_{i=1}^{n} \ln[1-h_{i}(x,\lambda ,q)^{r} ].
\end{eqnarray}
Now differentiating $ \ell(\theta) $ with respect to $ \alpha, \lambda,$ and $ q $ we get,
\begin{eqnarray}
\nonumber
\frac{\partial \ell (\theta)}{\partial \alpha}&=&\frac{n}{\lambda} + \sum_{i=1}^{n} \ln [1-h_{i}(x,\lambda ,q)^{r} ]\\
\nonumber
\frac{\partial \ell (\theta)}{\partial \lambda}&=&\frac{n}{\lambda} - (r-1)(1-q) \sum_{i=1}^{n} \frac{x_{i}}{h_{i}(x,\lambda ,q)} + r (\alpha -1) (1-q)\sum_{i=1}^{n} x_{i} \frac{h_{i}(x,\lambda ,q)^{r-1}}{1-h_{i}(x,\lambda ,q)^{r}} \\
\nonumber
\frac{\partial \ell (\theta)}{\partial q}&=&-\frac{n}{2-q}+\frac{r}{(2-q)(1-q)}\sum_{i=1}^{n}\ln(1-h_{i}(x,\lambda ,q))+(r-1)\sum_{i=1}^{n}\frac{\lambda x_{i} }{h_{i}(x,\lambda ,q)}\\
\nonumber
& &-(\alpha -1)\frac{r}{(2-q)(1-q)}\sum_{i=1}^{n}\frac{h_{i}(x,\lambda ,q)^{r} \ln(h_{i}(x,\lambda ,q))}{1-h_{i}(x,\lambda ,q)^{r}}\\
\nonumber
& &-(\alpha -1)r\sum_{i=1}^{n} \frac{\lambda x_{i}h_{i}(x,\lambda ,q)^{r-1} }{1-h_{i}(x,\lambda ,q)^{r}}
\end{eqnarray} \\
Now setting
$\frac{\partial \ell (\theta)}{\partial \alpha}=0 $,
$ \frac{\partial \ell (\theta)}{\partial \lambda}=0 $,
and $ \frac{\partial \ell (\theta)}{\partial q}=0 $,
and solving these system of equations simultaneously, we get the maximum likelihood estimate
$\hat{\theta}=(\hat{\alpha},\hat{\lambda},\hat{q})^T$ of $\theta=(\alpha,\lambda,q)^T$. For solving these non-linear equations we can use any iteration method such as Newton-Raphson technique.

\vskip.2cm
\subsection{Bayesian Estimation}
\vskip.2cm

In this section our focus is  to obtain the estimates of shape parameter of TLGqE distribution using  Bayesian
paradigm techniques by normal approximation.
Large sample Bayesian methods are primarily based on normal approximation
to the posterior distribution of $\theta$. As the sample size $n$ increases, the posterior distribution approaches normality under certain regularity conditions and
hence can be well approximated by an appropriate normal distribution if $n$ is
sufficiently large. When $n$ is large, the posterior distribution becomes highly
concentrated in a small neighborhood of the posterior mode, $\hat{\theta}$, for more details see Ghosh et al. (2006).
If the posterior distribution $f(\theta|x)$ is unimodal and roughly symmetric, it is convenient to
approximate it by a normal distribution centered at the mode, and  the logarithm of the posterior is approximated by a
quadratic function, yielding the approximation
$$f(\theta|x)\sim N(\hat{\theta}, [I(\hat{\theta})]^{-1} ), [I(\hat{\theta})]= -\frac{\partial^2 \ln
f(\theta|x) }{\partial \theta^2}$$
 If the mode, $\hat{\theta}$
is in the interior parameter space, then $I(\hat{\theta})$
is positive; if $\hat{\theta}$
is a vector parameter, then $I(\hat{\theta})$
is a
matrix. The estimation of shape parameter of TL
distribution using various Bayesian approximation techniques like normal approximation,
Tierney and Kadane (T-K) Approximation are given by Sultan and Ahmad (2015).

In our study the normal approximations of Topp-Leone distribution under different priors is obtained as under:
The likelihood function of (\ref{eq:lf}) for a sample of size $n$ is given as
\begin{equation} \label{eq:34}
L(\underline{x}|\alpha)\propto (\alpha)^n
{\rm e}^{-\alpha\Sigma_{i=1}^n\ln\left(2 G(x_i)-G(x_i)^2\right)^{-1}},
\end{equation}
where $G(x_i)=\left\{1-[1-(1-q) \lambda x_i]^{2( \frac{2-q}{1-q})}\right\}^{\alpha}$.
 Under uniform prior $g(\alpha)\propto 1$, the posterior distribution for $\alpha$ is given as
$$f(\alpha|x) \propto (\alpha)^n
{\rm e}^{-\alpha S},S=\Sigma_{i=1}^n\ln\left(2 G(x_i)-G(x_i)^2\right)^{-1}$$
and
$$\ln f(\alpha|x) = \ln K+ n \ln \alpha -\alpha S,$$ where  K is a constant. Then
$$\frac{\partial \ln
f(\alpha|x) }{\partial \alpha}=\frac{n}{\alpha}-S.$$ Hence the posterior mode is obtained as $\hat{\alpha}= \frac{n}{S}$ and
$I(\hat{\alpha})= \frac{S^2}{n}$.
Thus, the posterior distribution can be approximated as $$f(\alpha|x )\sim N \left(\frac{n}{S}, \frac{n}{S^2}\right).$$
Under extension of Jeffrey's prior $g(\alpha)\propto \left( \frac{1}{\alpha}\right)^m, m \in \mathcal{R^{+}}$ then the posterior distribution can be approximated as $$f(\alpha|x )\sim N \left(\frac{n-m}{S}, \frac{n-m}{S^2}\right).$$
Under gamma prior $g(\alpha)\propto {\rm e}^{-a \alpha}\alpha^{p-1}, a, p, \alpha >0$ then the posterior distribution can be approximated as
   $$f(\alpha|x )\sim N \left(\frac{n+p-1}{S+a}, \frac{n+p-1}{(S+a)^2}\right).$$

\vskip.2cm
\subsection{Random Variate Generation}
\vskip.2cm
\noindent
By using inversion method, we can generate a randon variate from TLqE distribution. We have already seen  that the relationship between a random variable $X$, having TLqE distribution, and a random variable $T$, having the TL distribution, is
\begin{eqnarray}\label{eq:31}
\nonumber
X&=&G^{-1}(t)\\
&=&\frac{1-(1-t)^{\frac{1-q}{2-q}}}{(1-q) \lambda}
\end{eqnarray}\\
where $G^{-1}(\cdotp)$ is related to inversion of the $q$-exponential cdf. The quantile function of the TL distribution is
\begin{equation} \label{eq:32}
t=1-\sqrt{1-u^{\frac{1}{\alpha}}},
\end{equation}
where $u$ is picked from the uniform distribution over $(0,1)$. Then the quantile function of TLqE distribution is obtained by substituting equation (\ref{eq:32}) into equation (\ref{eq:31}),
\begin{eqnarray*}
	X&=&\frac{1-\left(\sqrt{1-u^{\frac{1}{\alpha}}}\right)^{(\frac{1-q}{2-q})}}{(1-q)\lambda}
\end{eqnarray*}
For example, if $u$=0.7235, 0.9690, 0.5374, 0.8221, 0.1961, then $X \sim $ TLqE(0.3,1.5,1.2) can be generated respectively as
\begin{eqnarray*}
	x_{1}&=&\frac{1-\left(\sqrt{1-(0.7235)^{\frac{1}{0.3}}}\right)^{(\frac{1-1.2}{2-1.2})}}{(1-1.2)1.5}=0.1777\\ \\
	x_{2}&=&\frac{1-\left(\sqrt{1-(0.9690)^{\frac{1}{0.3}}}\right)^{(\frac{1-1.2}{2-1.2})}}{(1-1.2)1.2}=1.1137\\ \\
	x_{3}&=&\frac{1-\left(\sqrt{1-(0.5374)^{\frac{1}{0.3}}}\right)^{(\frac{1-1.2}{2-1.2})}}{(1-1.2)1.2}=0.0567\\ \\
	x_{4}&=&\frac{1-\left(\sqrt{1-(0.8221)^{\frac{1}{0.3}}}\right)^{(\frac{1-1.2}{2-1.2})}}{(1-1.2)1.2}=0.3201\\ \\
	x_{5}&=&\frac{1-\left(\sqrt{1-(0.1961)^{\frac{1}{0.3}}}\right)^{(\frac{1-1.2}{2-1.2})}}{(1-1.2)1.2}=0.0017\\ \\
\end{eqnarray*}
Thus using this technique we can simulate random variates for any values of the parameters.

\vskip.2cm
\subsection{Application}
\vskip.2cm
\noindent In this section, we consider a real data set and try to find the distribution that fits better to the data among the TLqE distrubution and TLE distribution. For the purpose of model selection, we use the Akaike Information Criterion (AIC), and the Bayesian Information criterion (BIC). The real data set that we consider, is the ball bearing data, which says the number of revolutions before failure for ball bearing (Crowder et al.,1994).

The data is
33.00, 68.64, 173.40, 41.52, 42.12, 68.64, 68.88, 45.60, 48.48,
84.12, 93.12,  98.64, 105.12, 105.84, 51.84, 51.96, 54.12, 17.88,
55.56, 127.92, 128.04, 67.80,  67.80, 28.92.
 We provide the values of estimated parameters, AIC and BIC values for ball bearing data in the table given below.

\begin{center}
\begin{tabular}{|c|ccc|cc|}
\hline
The model & \multicolumn{3}{c|}{Estimate}&AIC&BIC\\
&$\alpha$&$\lambda$&$ q $ &&\\

\hline
TLE&5.2827&0.0161&&229.9559&119.2489\\
TLqE&2.4037&0.0118&1.1710&225.1589&118.9860\\
\hline
\end{tabular}
\end{center}
The AIC and BIC values of the TLqE distribution is the smallest ones among the three considered distributions. As expected TLqE is more appropriate for this data than the TLE distribution because of its shape of hazard function.


\end{document}